\documentclass[12pt]{article}
\usepackage{amssymb}
\usepackage{amsmath}
\usepackage{indentfirst,latexsym}
\usepackage{mathrsfs}
\usepackage{graphicx}
\usepackage{fancyhdr}

\date{}
\begin{document}

\title{Restricted involutions and Motzkin paths}
\author{Marilena Barnabei, Flavio Bonetti, and Matteo Silimbani \thanks{
Dipartimento di Matematica - Universit\`a di Bologna}} \maketitle

\noindent {\bf Abstract.} We show how a bijection due to Biane
between involutions and labelled Motzkin paths yields bijections
between Motzkin paths and two families of restricted involutions
that are counted by Motzkin numbers, namely, involutions avoiding
$4321$ and $3412$. As a consequence, we derive characterizations
of Motzkin paths corresponding to involutions avoiding either
$4321$ or $3412$ together with any pattern of length 3.
Furthermore, we exploit the described bijection to study some
notable subsets of the set of restricted involutions, namely,
fixed point free and centrosymmetric restricted involutions.
\newline

\newline


\section{Introduction}

\newtheorem{yama}{Proposition}

\noindent A permutation $\sigma\in S_n$ avoids the pattern
$\tau\in S_k$ if $\sigma$ does not contain a subsequence
order-isomorphic to $\tau$. Permutations with forbidden patterns
have been intensively studied in recent years for their connection
with many problems arising in both computer science and
combinatorics, such as Schubert varieties, Kazdan-Lusztig
polynomials, Chebyshev polynomials and various sorting algorithms
(see \cite{gtouf}, \cite{kratt}, \cite{wtouf}, \cite{simi}, and
references there in). Many classical sequences, such as Catalan
numbers, Motzkin numbers, central binomial coefficients and
Fibonacci numbers, occur when the cardinalities of sets of pattern
avoiding permutations are computed.

\noindent In particular, Motzkin numbers arise mainly when
involutions avoiding certain patterns of length $4$ are studied.
More precisely, involutions on $n$ objects avoiding the patterns
$3412$, $4321$, $2143$, and $1234$ enumerated by Motzkin numbers
(see \cite{gui}, \cite{ppg}, and \cite{reg}).\newline

\noindent On the other hand, in the literature there are many
bijections between permutations and Motzkin paths provided with
some kind of labelling (see \cite{bia}, \cite{bmv}, \cite{fz}, and
\cite{fv}). In this paper, we consider the bijection defined by
Biane in \cite{bia}, restricted to the set of involutions. The
bijection maps an involution into a Motzkin path whose down steps
are labelled with an integer that does not exceed its height,
while the other steps are unlabelled.

\noindent This bijection reveals to be an effective tool for
characterizing several sets of restricted involutions, and allows
us both to recover already known results in a unified framework,
and to obtain some new characterizations.

\noindent In particular, we prove that an involution $\tau$ avoids
$4321$ if and only if the label of any down step in the corresponding Motzkin
path is 1. Similarly, an involution $\tau$
avoids $3412$ if and only if the label of any down step in the corresponding
Motzkin path equals its height.

\noindent These two results allow us, on the one hand, to regain
the two standard bijections between (unlabelled) Motzkin paths and
involutions avoiding either $4321$ or $3412$ (see \cite{gui}),
and, on the other hand, to characterize those labelled Motzkin
paths that correspond to involutions avoiding one of these two
patterns together with an assigned pattern of length $3$, hence
obtaining some new enumerative results.\noindent

\noindent In the last two sections, as an example of the efficacy
of the present approach, we study pattern avoidance on sets of
involutions corresponding to two notable subsets of Motzkin paths,
namely, Dyck paths and Motzkin paths that are symmetric with
respect to a vertical line.

\section{Preliminary notions}

\noindent Given a permutation $\sigma\in S_n$, one can partition
the set $\{1,2,\ldots,n\}$ into intervals $A_1,\ldots,A_t$ (i.e.
$A_i=\{k,k+1,\ldots,k+h\}$) such that $\sigma(A_j)=A_j$ for every
$j$. The restrictions of $\sigma $ to the intervals in the finest
of these decompositions are called \emph{connected components} of
$\sigma$. A permutation $\sigma$ with a single connected component
is called \emph{connected}.

\noindent For example, the permutation $$\sigma=2\ 4\ 5\ 3\ 1\ 7\
6$$ has the two connected components $\sigma_1=2\ 4\ 5\ 3\ 1$ and
$\sigma_2=7\ 6$, while the permutation
$$\rho=2\ 7\ 6\ 1\ 3\ 5\ 4$$
is connected. \newline

\noindent Let $\psi\in S_n$ be defined by $\psi(i)=n+1-i$. If
$\sigma\in S_n$, the \emph{reverse} of $\sigma$ is the permutation
$\sigma_r=\sigma\psi$.  Similarly, the \emph{complement} of
$\sigma$ is the permutation $\sigma_c=\psi\sigma$ and the
\emph{reverse-complement} of $\sigma$ is the permutation
$\sigma_{rc}=\psi\sigma\psi$.\\

\noindent Let $\sigma\in S_n$ and $\pi\in S_k$, $k\leq n$, be two
permutations. We say that $\sigma$ \emph{contains the pattern}
$\pi$ if there exists a subsequence
$(\sigma(i_1),\sigma(i_2),\ldots,\sigma(i_k))$ with $1\leq
i_1<i_2<\cdots<i_k\leq n$ that is order-isomorphic to $\pi$. We
say that $\sigma$ \emph{avoids} the pattern $\pi$ if $\sigma$ does
not contain $\pi$. Given a set of permutations $A\subseteq S_n$
and a set of patterns $\pi_1,\ldots,\pi_k$, we denote by
$A(\pi_1,\ldots,\pi_k)$ the set of elements in $A$ that avoid
$\pi_i$ for every $1\leq i\leq k$.

\noindent We say that two patterns $\pi_1$ and $\pi_2$ are
\emph{equidistributed} over a set of permutations $A$ if
$|A(\pi_1)|=|A(\pi_2)|$ and that they are \emph{equivalent} over
$A$ if
$A(\pi_1)=A(\pi_2)$.\\

\noindent Recall that a permutation $\tau\in S_n$ is an involution
if and only if $\tau^{-1}=\tau$. Equivalently, $\tau$ is an
involution if and only if its cycle structure has no cycle of
length longer than two.

\noindent Denote by $I_n$ the set of involutions in the symmetric
group $S_n$. We point out that the bijection that associates a
permutation $\sigma$ to its reverse-complement $\sigma_{rc}$ maps
$I_n$ into itself. This implies that a pattern $\pi\in S_k$ and
its reverse-complement $\pi_{rc}$ are equidistributed over $I_n$.
In fact, a permutation $\sigma$ contains $\pi$ if and only if
$\sigma_{rc}$ contains $\pi_{rc}$.

\noindent Moreover, any pattern $\pi$ is equivalent over $I_n$ to
its inverse $\pi^{-1}$, since, in general, for every permutation
$\sigma\in S_n$, $\sigma$ contains $\pi$ if and only if
$\sigma^{-1}$ contains $\pi^{-1}$.

%\noindent These observations reduce the number of cases that we
%need to consider in the present instance.\\

\section{Labelled Motzkin paths and involutions}

\noindent A \emph{Motzkin path} of length $n$ is a lattice path
starting at $(0,0)$, ending at $(n,0)$, and never going below the
$x$-axis, consisting of up steps $U=(1,1)$, horizontal steps
$H=(1,0)$, and down steps $D=(1,-1)$. A \emph{Dyck path} is a
Motzkin path that does not contain horizontal steps.

\noindent The set of Motzkin paths of length $n$ will be denoted
by $\mathscr{M}_n$, while the set of Dyck paths of length $2n$
will be denoted by $\mathscr{D}_n$. It is well known that the
cardinality of $\mathscr{M}_n$ is the $n$-th \emph{Motzkin number}
$M_n$ (sequence A001006 in \cite{oeis}), and that the cardinality
of the set $\mathscr{D}_n$ is
the $n$-th Catalan number $C_n$ (sequence A000108 in \cite{oeis}).\\

\noindent We define the \emph{height of a step} of a Motzkin path
to be the larger $y$ coordinate of the step, and the \emph{height
of a path} to be the largest height of its steps.\\

\noindent An \emph{irreducible} (or \emph{elevated}) Motzkin path
is a Motzkin path that does not touch the $x$-axis except for the
origin and the final destination. An \emph{irreducible component}
of a Motzkin path $M$ is a
maximal irreducible Motzkin subpath of $M$.\\

\noindent A \emph{labelling} of a Motzkin path $M$ will be a map
that associates with every down step $D$ at height $h$ an integer
$\lambda(D)$ such that $1\leq\lambda(D)\leq h$. A \emph{labelled
Motzkin path} is a pair $(M,\lambda)$, where $M$ is a Motzkin path
and $\lambda$ a labelling of $M$.

\noindent We will denote by $\upsilon$ the \emph{unitary
labelling}, namely, the labelling that assigns the label $1$ to
every down step $D$ of $M$. Similarly, the \emph{maximal
labelling}
$\mu$ will be defined by $\mu(D)=h(D)$, where $h(D)$ is the height of $D$.\\

\noindent We now define a bijection $\Phi$ between the set $I_n$
of involutions of the symmetric group $S_n$ and the set of
labelled Motzkin paths of length $n$, that is essentially a
restriction to the set of involution of the bijection appearing in
\cite{bia}.

\begin{figure}[ht]
\begin{center}
\includegraphics[bb=101 638 374 735,width=.5\textwidth]{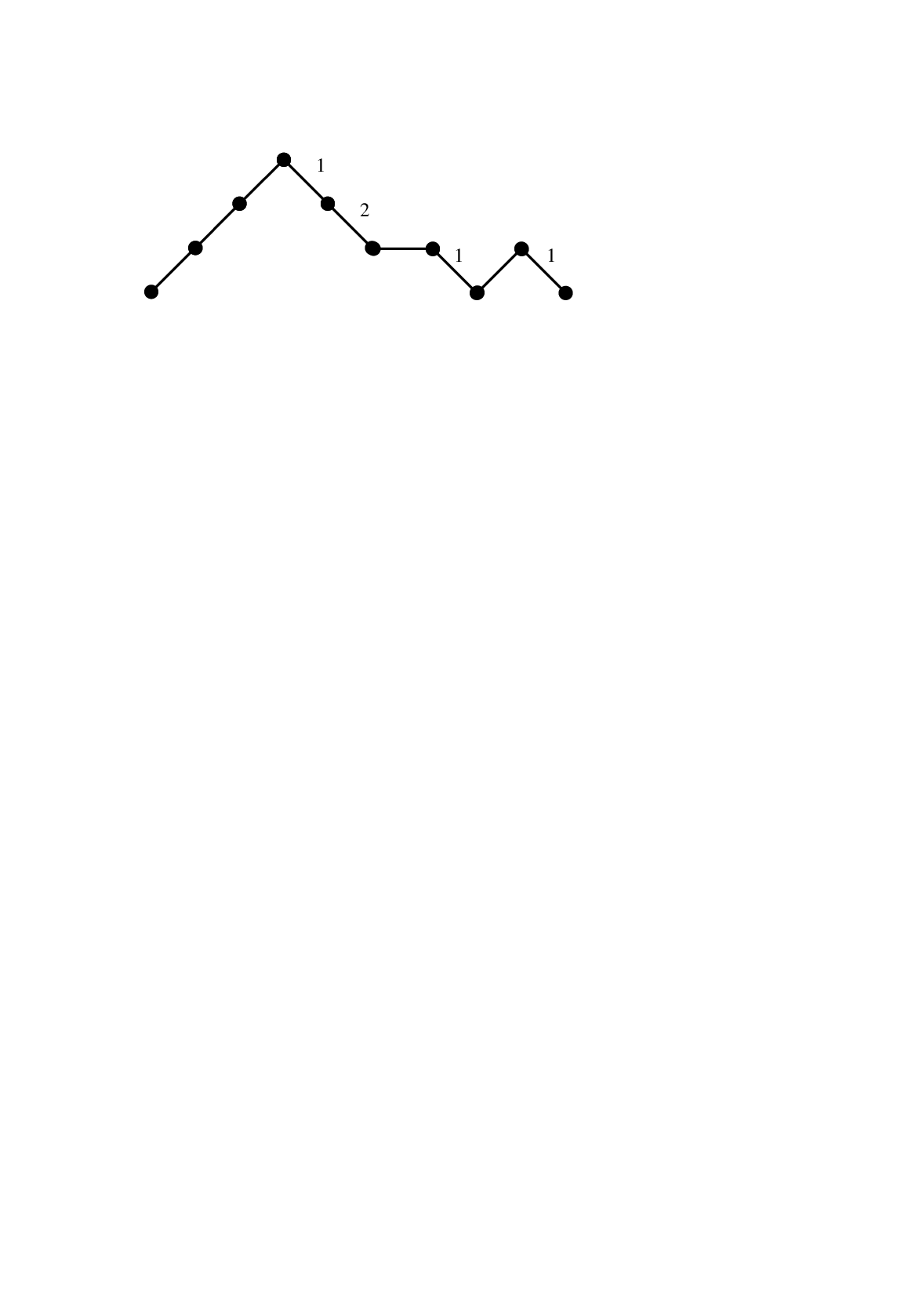} \caption{The labelled Motzkin path $\Phi(\tau)$, with $\tau=4\ 7\ 5\ 1\ 3\ 6\
2\ 9\ 8$}\label{labmot}
\end{center}
\end{figure}

\noindent Consider an involution $\tau$ and determine the set
exc$(\tau)=\{i|\tau(i)>i\}$ of its \emph{excedances}. Starting
from the empty path and from the list $A_{\tau}=(a_1\ a_2\ \ldots\
a_s)$ consisting of the elements of exc$(\tau)$ written in
increasing order, we construct the labelled Motzkin path
$\Phi(\tau)$ by adding a step for every integer $1\leq i\leq n$ as
follows:

\begin{itemize}
\item[-] if $\tau(i)=i$ we add a horizontal step at $i$-th position;
\item[-] if $\tau(i)>i$ we add an up step at $i$-th position;
\item[-] if $\tau(i)<i$ we add a down step at $i$-th position. The
integer $\tau(i)$ appears in the list $A_{\tau}$ at position $k$,
say. We assign the label $k$ to the new down step and remove
$\tau(i)$ from $A_{\tau}$.
\end{itemize}
For example, the image $\Phi(\tau)$ of the involution $\tau=4\ 7\
5\ 1\ 3\ 6\ 2\ 9\ 8$ is the labelled Motzkin path in Figure
\ref{labmot}.

\noindent Note that an involution $\tau$ is connected if and only
if its image is an irreducible Motzkin path.\\

\noindent There is a simple relationship between the Motzkin paths
associated with an involution $\tau$ and its reverse-complement
$\tau_{rc}$:

\newtheorem{evacu}[yama]{Theorem}
\begin{evacu}\label{modo}
Let $(M,\lambda)$ and $(M',\lambda')$ the labelled Motzkin paths
of length $n$ associated by $\Phi$ with $\tau$ and $\tau_{rc}$,
respectively. Then, the Motzkin paths $M$ and $M'$ have the
following symmetry (with respect to the vertical line
$x=\frac{n}{2}$):\begin{itemize}
\item[i.] if $M$ has a horizontal step at the $k$-th position,
then $M'$ has a horizontal step at the $(n+1-k)$-th position;
\item[ii.] if $M$ has an up step at the $k$-th position, then $M'$ has a
down step at the $(n+1-k)$-th position;
\item[iii.] if $M$  has a down step at the $k$-th position, then $M'$ has
an up step at the $(n+1-k)$-th position.
\end{itemize}
\end{evacu}

\noindent \emph{Proof} Case $i.$ is easily verified. For case
$ii.$, suppose that $M$ has an up step at the $k$-th position.
This implies that $j=\tau(k)>k$. Then we must have
$\tau_{rc}(n+1-k)=n+1-j<n+1-k$, and hence the $n+1-k$-th step in
$M'$ is a down step. The last case follows analogously.
\begin{flushright}
$\diamond$
\end{flushright}

\noindent Motzkin numbers appear in the enumeration of involutions
avoiding several patterns of length $4$. More precisely, for
$\tau$ equal to $1234$, $1243$, $3412$, $3214$, $1432$, $4321$ or
$2143$, $|I_n(\tau)|=M_n$ (see \cite{chen}, \cite{egg},
\cite{gui}, \cite{ppg}, and \cite{reg}). Hence, it is natural to
search for a bijection between the sets $I_n(\tau)$ and
$\mathscr{M}_n$ for each one of these patterns. Some of these
bijections appear in the literature (see e.g. \cite{chen} and
\cite{egg}). We will show that the described correspondence
between $I_n$ and the set of labelled Motzkin paths of length $n$
provides in a natural way such a bijection in the cases
$\tau=3412$ and $\tau=4321$.

\section{$4321$-avoiding involutions}

\newtheorem{qtdu}[yama]{Theorem}
\begin{qtdu}\label{nlsn}
Let $\tau$ be an involution and $\Phi(\tau)=(M,\lambda)$. Then,
$\tau$ avoids the pattern $4321$ if and only if $\lambda$ is the
unitary labelling.
\end{qtdu}

\noindent \emph{Proof} Suppose that there exists a down step at
position $i$ in $M$ with a label greater than $1$. Clearly, such a
step can not be the last one in $M$. Let $j$ be the least position
greater than $i$ corresponding to a down step whose label is $1$
(in the worst case, $j=n$). Then, $\tau$ contains the subsequence
$$j\ i\ \tau(i)\ \tau(j)$$
which is of type $4321$.

\noindent Conversely, suppose that each label in $M$ equals $1$.
In this case, $\tau$ is obtained by interlacing the three
following increasing subsequences:

\begin{itemize}
\item[-] the sequence of fixed points;
\item[-] the sequence of excendances;
\item[-] the sequence of deficiencies, namely, the sequence of the images of excedances of
$\tau$,
\end{itemize}
and, hence, $\tau$ can not contain any decreasing subsequence of
length $4$.
\begin{flushright}
$\diamond$
\end{flushright}

\begin{figure}[ht]
\begin{center}
\includegraphics[bb=104 582 468 714,width=.7\textwidth]{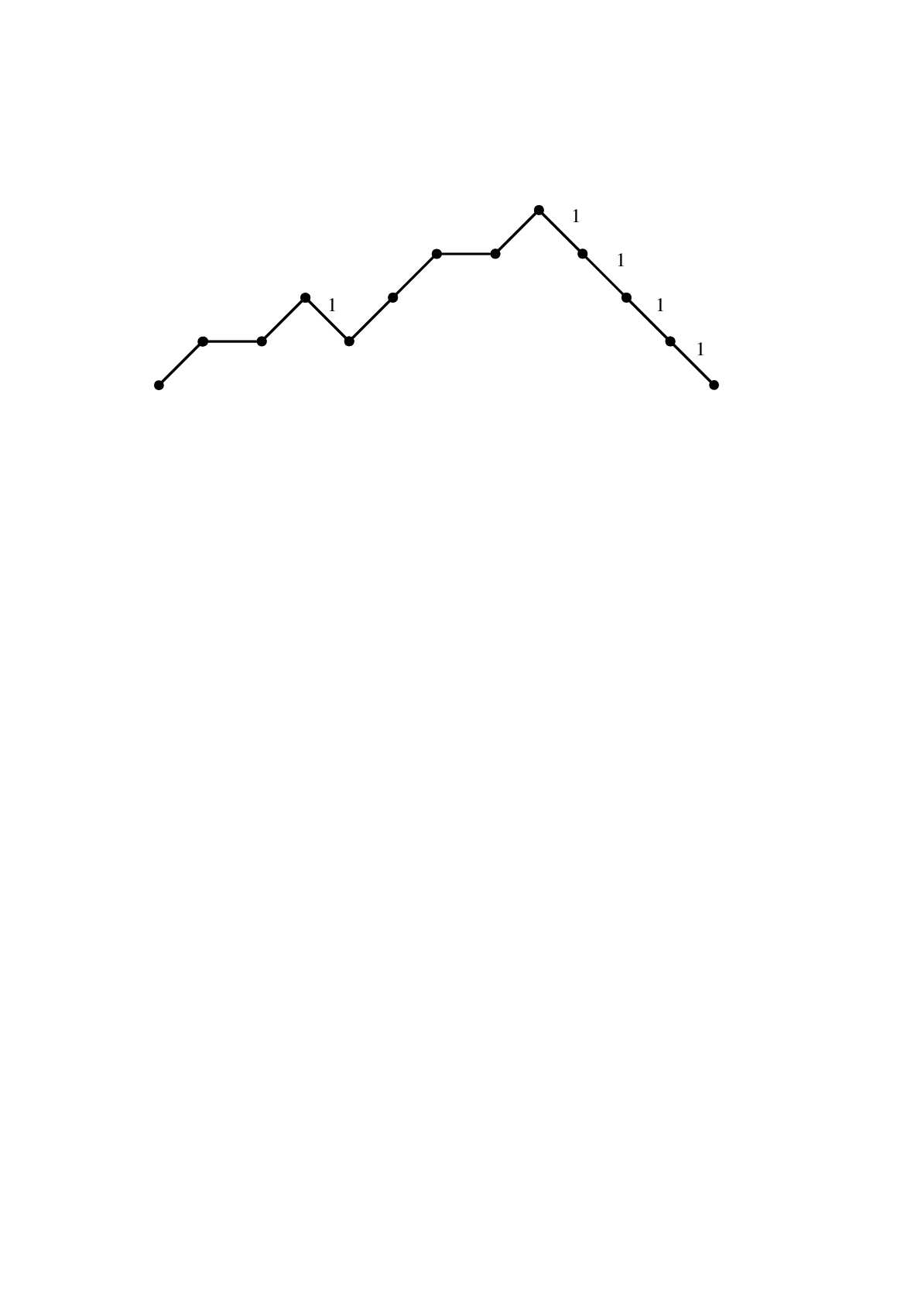}
\caption{The labelled Motzkin path associated with the $4321$
avoiding involution $\tau=4\ 2\ 9\ 1\ 10\ 11\ 7\ 12\ 3\ 5\ 6\ 8
$}\label{fors}
\end{center}
\end{figure}

\noindent We remark that also the involutions in $I_n(321)$ and in
$I_n(312)$ have a simple characterization in term of labelled
Motzkin paths:

\newtheorem{tdu}[yama]{Proposition}
\begin{tdu}\label{replace}
Let $\tau$ be an involution in $I_n$ and $\Phi(\tau)=(M,\lambda)$.
Then, $\tau$ avoids $321$ if and only if $\lambda=\upsilon$ and
all horizontal steps in $M$ are at height $0$.
\end{tdu}

\noindent \emph{Proof} \noindent Suppose that $M$ contains a
horizontal step $H$ at position $k$ and height $h>0$. In this
case, consider the up step $U$ preceding $H$ and closest to $H$.
Denote by $j$ the position of $U$ in $M$. The integer $j$
corresponds to the position of an excedance of $\tau$, so the
involution $\tau$ contains the $321$-sequence $\tau(j)\ k\ j$.

\noindent Suppose now that the labelling $\lambda$ is not unitary.
As seen above, in this case $\tau$ contains a $4321$-subsequence,
and hence contains also the pattern $321$.

\noindent Conversely, suppose that $M$ has no horizontal steps at
a height greater than $0$ and that $\lambda$ is unitary. Then, the
irreducible components of $M$ are either horizontal steps at
height $0$ or Dyck paths. Recall that the irreducible components
of $M$ correspond to the connected components of $\tau$. The
connected components corresponding to horizontal steps avoid
trivially the pattern $321$. On the other hand, each component
corresponding to a Dyck path (with the unitary labelling) can be
obtained by interlacing the increasing sequences of its excedances
and deficiencies. This completes the proof.
\begin{flushright}
$\diamond$
\end{flushright}

\begin{figure}[ht]
\begin{center}
\includegraphics[bb=158 636 403 710,width=.5\textwidth]{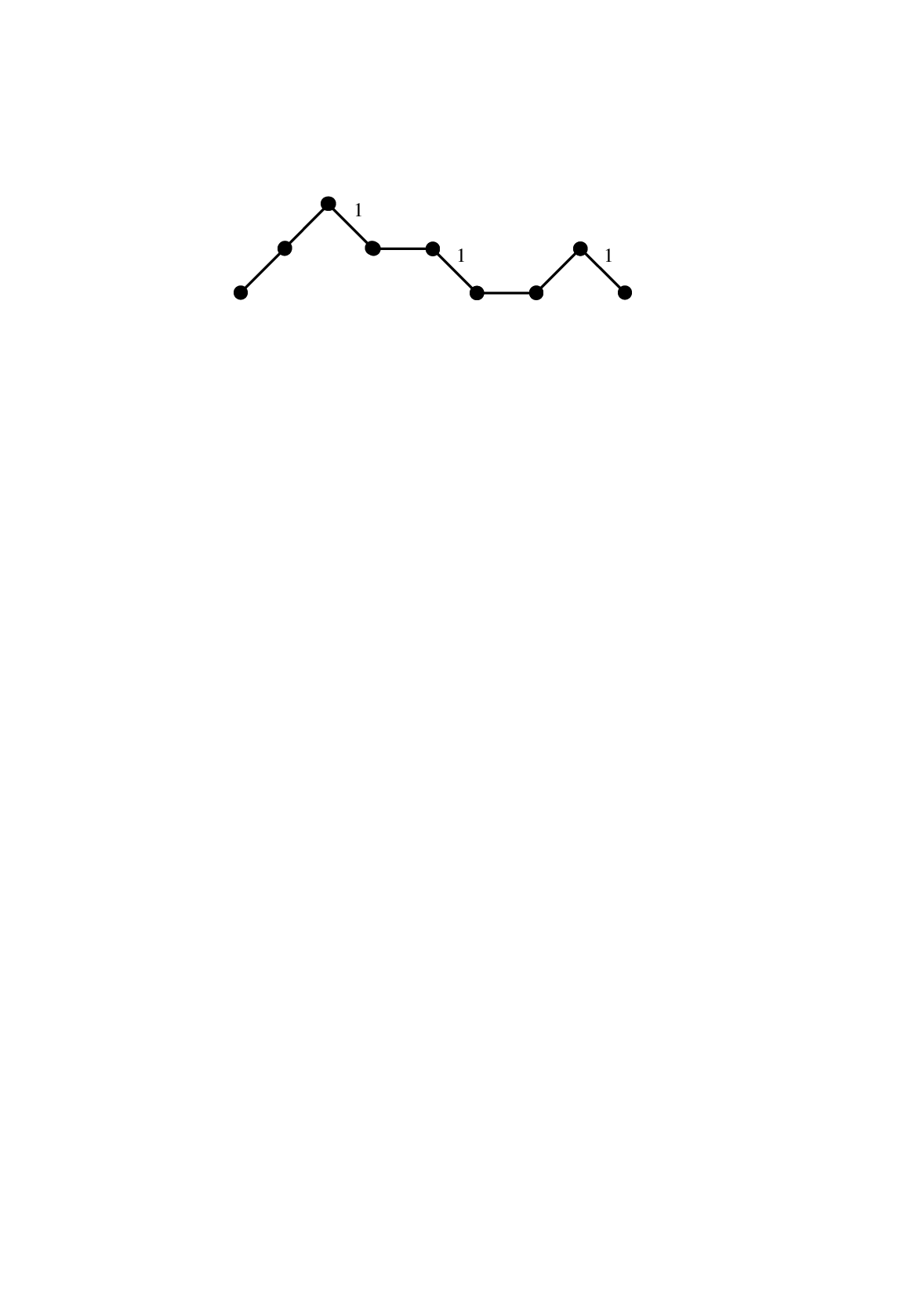} \caption{The labelled Motzkin path corresponding to $\tau=3\ {\bf 5}\ 1\ {\bf 4}\ {\bf 2}\
6\ 8\ 7$, that contains $321$.}
\end{center}
\end{figure}

\newtheorem{tud}[yama]{Proposition}
\begin{tud}\label{fumetti}
Let $\tau$ be an involution in $I_n$ and $(M,\lambda)$ the
associated labelled Motzkin path. Then, $\tau$ avoids $312$ if and
only if $\lambda=\mu$ and the irreducible components of $M$ are of
 type either $U^aHD^a$ or $U^aD^a$.
\end{tud}

\noindent \emph{Proof} Remark that an involution $\tau$ avoids
$312$ if and only if there exist integers $k_1<k_2<\cdots<k_p$
such that $\tau$ can be written in one-line notation as follows
$$k_1\ k_1-1\ \cdots\ 1\ k_2\ k_2-1\ \cdots\ k_1+1\ \cdots\ k_p\ \cdots\ k_{p-1}+1.$$
In fact, suppose that $\tau(1)=k_1$. Hence $\tau(k_1)=1$, since
$\tau$ is an involution. As $\tau$ avoids $312$, we must have
$$\tau(2)=k_1-1,\ \tau(3)=k_1-2,\ \cdots\ \tau(k_1-1)=2.$$ Iterating this argument
we get the assertion.

\noindent This implies that involutions in $I_n(312)$ correspond
via the bijection $\Phi$ to Motzkin paths whose irreducible
components are of type either $U^aHD^a$ or $U^aD^a$, with the
maximal labelling.
\begin{flushright}
$\diamond$
\end{flushright}

\begin{figure}[ht]
\begin{center}
\includegraphics[bb=67 583 387 627,width=.7\textwidth]{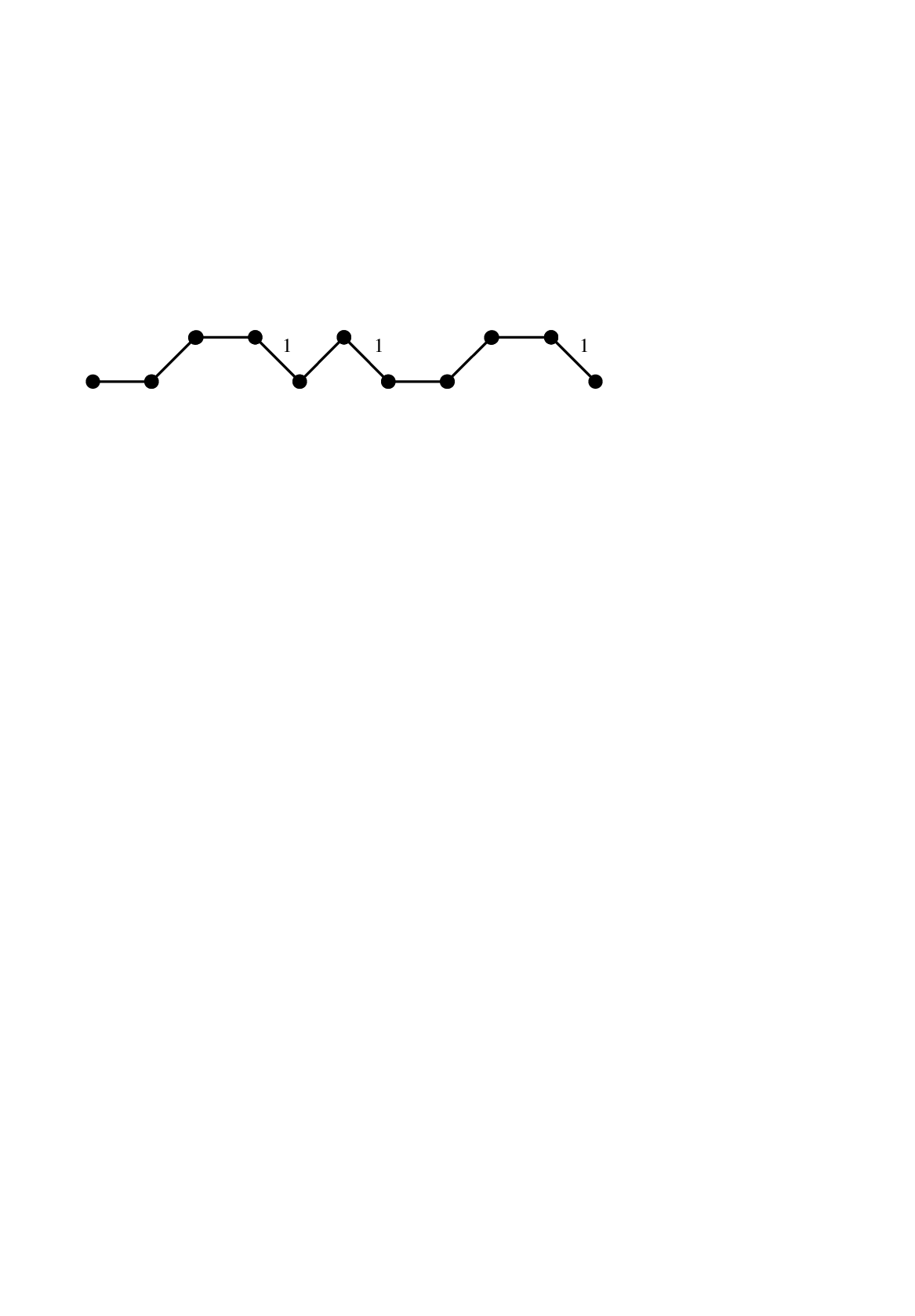} \caption{The labelled Motzkin path corresponding to $\tau=1\ 4\ 3\ 2\ 6\
5\ 7\ 10\ 9\ 8$, that avoids $312$.}
\end{center}
\end{figure}

\noindent We now want to characterize labelled Motzkin paths
corrosponding to involutions avoiding $4321$ and a pattern
$\sigma$ of length $3$. We do not consider the case $\sigma=123$
since, for every $n>6$, $I_n(4321,123)=\emptyset$.

\newtheorem{multiple}[yama]{Theorem}
\begin{multiple}\label{attrezz}
Let $\tau$ be an involution in $I_n(4321)$ associated with the
labelled Motzkin path $(M,\upsilon)$. Then:
\begin{itemize}
\item[i.] $\tau$ avoids $132$ if and only if $M$ does not contain
any subpath among $HU$, $DU$ and $DHD$. As a consequence, $\tau$
avoids $213$ if and only if $M$ does not contain any subpath among
$DH$, $DU$ and $UHU$;
\item[ii.] $\tau$ avoids $321$ if and only if
 all horizontal steps in $M$ have height $0$;
\item[iii.] $\tau$ avoids $312$, and hence $231$, if and only if
the irreducible components of $M$ are either $H$, or $UHD$, or
$UD$.
\end{itemize}
\end{multiple}

\noindent \emph{Proof}
\begin{itemize}
\item[i.] Consider the pattern $132$. Let $(M,\upsilon)$ be a Motzkin path with the unitary labelling.
 If $M$ contains one subpath among
$HU$, $DU$ and $DHD$, then the corresponding involution $\tau$
must contain $132$. In fact:
\begin{itemize}
\item[1.] if the $k$-th step in $M$ is horizontal and the
$(k+1)$-th step is an up step, then $\tau(k)=k$ and
$\tau(k+1)=j>k+1$. Hence, $\tau$ contains the $132$-pattern $k\ j\
k+1$;
\item[2.] if the $k$-th step in $M$ is a down step and the
$(k+1)$-th step is an up step, then $\tau(k)=j<k$ and
$\tau(k+1)=m>k+1$. Hence, $\tau$ contains the $132$-pattern $j\ m\
k+1$;
\item[3.] if the $k$-th step in $M$ is a down step, the
$(k+1)$-th step is horizontal and the $(k+2)$-th step is a down
step, then $\tau(k)=j<k$, $\tau(k+1)=k+1$ and $\tau(k+2)=m$, with
$j<m<k$. Hence, $\tau$ contains the $132$-pattern $j\ k+1\ m$.
\end{itemize}
If $M$ does not contain any of those subpaths, then
$M=U^aH^bD^aH^c$ for suited non-negative integers $a,b,c$ such
that $2a+b+c=n$. A Motzkin path of this type corresponds to an
involution $\tau$ whose one-line notations is
$$
e_1\ \cdots\ e_a\ a+1\ \cdots\ a+b\ d_1\ \cdots\ d_b\ 2a+b+1\
\cdots\ n,$$ where $e_1\ldots e_a$ are the excedances of $\tau$
written in increasing order and $d_1\ldots d_b$ are the
deficiencies of $\tau$ in increasing order. In this case, $\tau$
clearly avoids $132$.

\noindent Remark now that the pattern $213$ is the
reverse-complement of $132$. As in the case of the whole set
$I_n$, the set $I_n(4321)$ is closed under reverse-complement. In
fact, let $\tau\in I_n(4321)$ and suppose that $\sigma=\tau_{rc}$
contains the pattern $4321$. Then, there exist integers $i<j<h<k$
such that $\sigma(i)>\sigma(j)>\sigma(h)>\sigma(k)$. Since
$\sigma(x)=n+1-\tau(n+1-x)$, the involution $\tau$ contains the
subsequence $n+1-k>n+1-h>n+1-j>n+1-i$, and this yields a
contradiction.  Then, if $\tau$ corresponds to the labelled
Motzkin path $(M,\upsilon)$, then $\tau_{rc}$ corresponds to the
path $M'$ obtained by reflecting $M$ over the line $x=\frac{n}{2}$
(see Theorem \ref{modo}), with the unitary labelling.

\begin{figure}[ht]
\begin{center}
\includegraphics[bb=42 535 535 628,width=.9\textwidth]{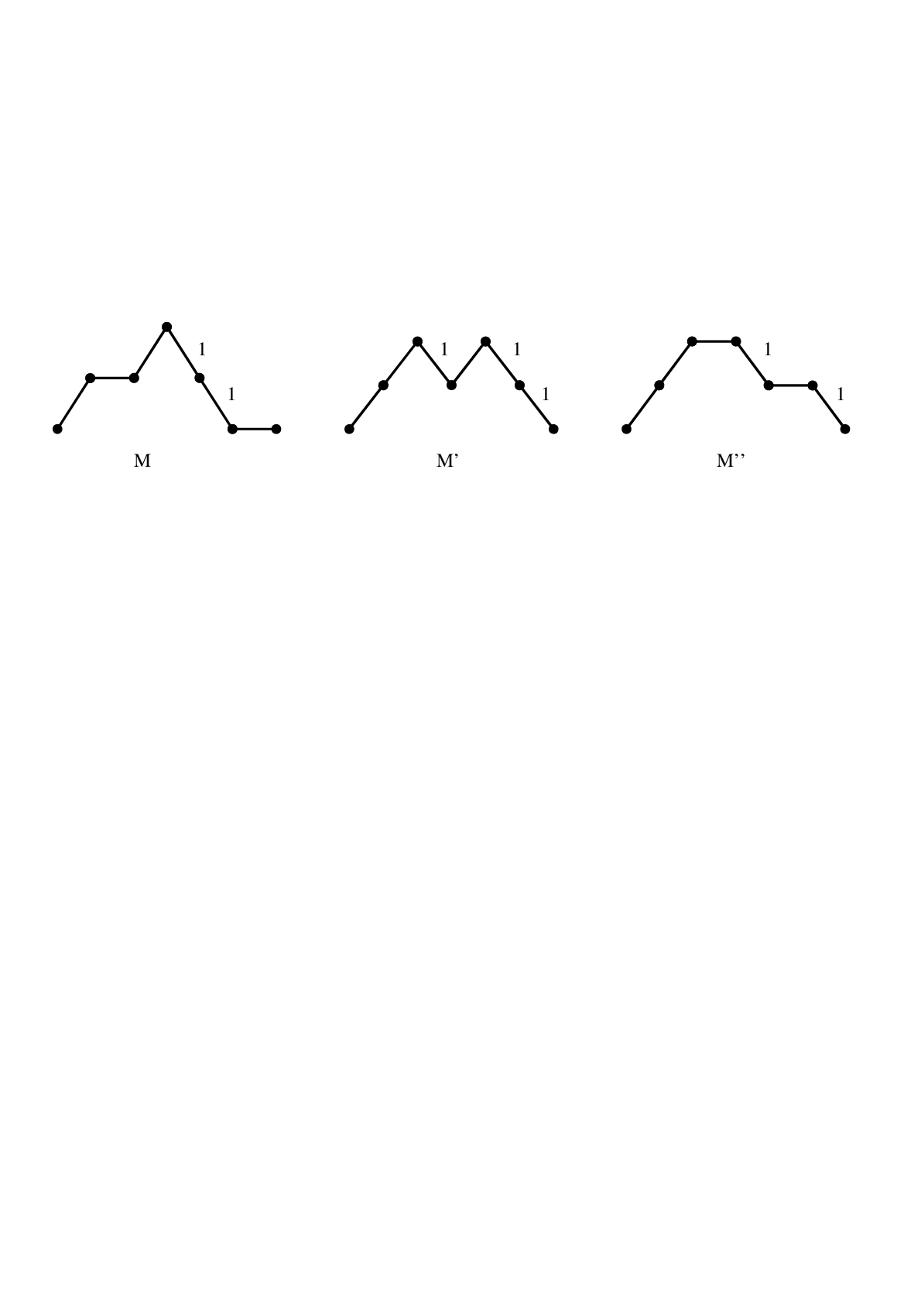} \caption{The labelled Motzkin paths $M$, $M'$, and $M''$ corresponding to
the involutions $\tau=4\ {\bf 2}\ {\bf 5}\ 1\ {\bf 3}\ 6$,
$\tau'=3\ 5\ {\bf 1}\ {\bf 6}\ 2\ {\bf 4}$, $\tau''=4\ 6\ 3\ {\bf
1}\ {\bf 5}\ {\bf 2}$, all containing $132$.}
\end{center}
\end{figure}

\item[ii.] The assertion follows immediately from Proposition
\ref{replace}, since an involution that avoids $321$ avoids also
$4321$.

\item[iii.] Recall that involutions in $I_n(312)$ correspond to Motzkin
paths whose irreducible components are either $B_a=U^aHD^a$ or
$P_a=U^aD^a$, with the maximal labelling (see Proposition
\ref{fumetti}). Such a path corresponds to an involution in
$I_n(4321)$ if and only if its height is $1$, namely, it is a
concatenation of subpaths of the kind $B_0=H$, $B_1=UHD$ and
$P_1=UD$.  This completes the proof.
\end{itemize}
\begin{flushright}
$\diamond$
\end{flushright}

\noindent The preceding theorem yields the following enumerative
results. Formulas $i.$ and $iii.$, up to our knowledge, are new.

\newtheorem{diecipunti}[yama]{Corollary}
\begin{diecipunti}
We have:
\item[i.] $|I_n(4321,132)|=|I_n(4321,213)|=1+\lfloor\frac{n}{2}\rfloor\lceil\frac{n}{2}\rceil$;
\item[ii.] $|I_n(4321,321)|={n\choose \lfloor\frac{n}{2}\rfloor}$;
\item[iii.] $I_n(4321,312)=I_n(4321,231)$ and
$|I_n(4321,312)|=t_{n+2}$,

\noindent where $t_n$ is the sequence of Tribonacci numbers
defined by $t_0=0$, $t_1=0$, $t_2=1$ and
$t_{n+3}=t_{n+2}+t_{n+1}+t_n$.
\end{diecipunti}

\noindent \emph{Proof}
\begin{itemize}
\item[i.] The cardinality of $I_n(4321,132)$
equals the number of Motzkin paths of the kind $M=U^aH^bD^aH^c$.
Hence, we have:
$$|I_n(4321,132)|=1+\sum_{h=1}^{\left\lfloor\frac{n}{2}\right\rfloor}n-2h+1=1+\lfloor\frac{n}{2}\rfloor\lceil\frac{n}{2}\rceil.$$
\item[ii.] Obviously $I_n(4321,321)=I_n(321)$. It is well known (see \cite{ss}) that the cardinality of
$I_n(321)$ equals the $n$-th central binomial coefficient.
\item[iii.] A Motzkin path associated with an involution in $I_n(4321,312)$ has irreducible components
 of the kind $H$, $B_1=UHD$ and $P_1=UD$. Hence, such a
Motzkin path $M\in\mathscr{M}_{n+3}$ can be obtained by adding
\begin{itemize}
\item[1.] an $H$ step at the end of a path of the same kind in
$\mathscr{M}_{n+2}$;
\item[2.] a subpath $P_1$ at the end of a path of the same kind in
$\mathscr{M}_{n+1}$;
\item[3.] a subpath $B_1$ at the end of a path of the same kind in
$\mathscr{M}_{n}$.
\end{itemize}
Moreover, $|I_1(4321,312)|=1=t_3$, $|I_2(4321,312)|=2=t_4$ and
$|I_3(4321,312)|=4=t_5$. This completes the proof.
\end{itemize}
\begin{flushright}
$\diamond$
\end{flushright}

\section{$3412$-avoiding involutions}

\noindent Also $3412$-avoiding involutions are enumerated by the
Motzkin numbers (see \cite{gui}). The labelled Motzkin paths
associated with an involution in $I_n(3412)$ can be characterized
as follows:

\newtheorem{tqud}[yama]{Theorem}
\begin{tqud}\label{consimm}
Let $\tau$ be an involution associated with the labelled Motzkin
path $(M,\lambda)$. Then, $\tau$ avoids the pattern $3412$ if and
only if $\lambda$ is the maximal labelling.
\end{tqud}

\noindent \emph{Proof} Let $\tau$ be the involution corresponding
to a given labelled Motzkin path and suppose that there exists a
down step $D$ at position $k$ whose label is not maximal. This
means that, at the $k$-th step of the procedure that creates the
Motzkin path $\Phi(\tau)$, the maximal element $a$ in the list
$A_{\tau}$ is not removed, and hence there exists an integer $j>k$
such that $\tau(j)=a$. This implies that the sequence
$\tau(1)\cdots\tau(n)$ contains the $3412$-subsequence $k\ j\
\tau(k)\ a$.

\noindent Conversely, suppose that the Motzkin path $\Phi(\tau)$
is maximally labelled. Suppose now that $\tau$ contains the
$3412$-subsequence $a\ b\ c\ d$.  Remark that, in this case, both
the subsequences of excedances and deficiencies of $\tau$ must be
decreasing. This implies that at least one among the integers $a$,
$b$ , $c$, and $d$ must be a fixed point of $\tau$. If $a$ or $b$
is fixed, then $c$ and $d$ are deficiencies in increasing order.
Similarly, if $c$ or $d$ is fixed, then $a$ and $b$ are
excendances. In both cases, we get a contradiction.
\begin{flushright}
$\diamond$
\end{flushright}

\begin{figure}[ht]
\begin{center}
\includegraphics[bb=104 582 468 714,width=.7\textwidth]{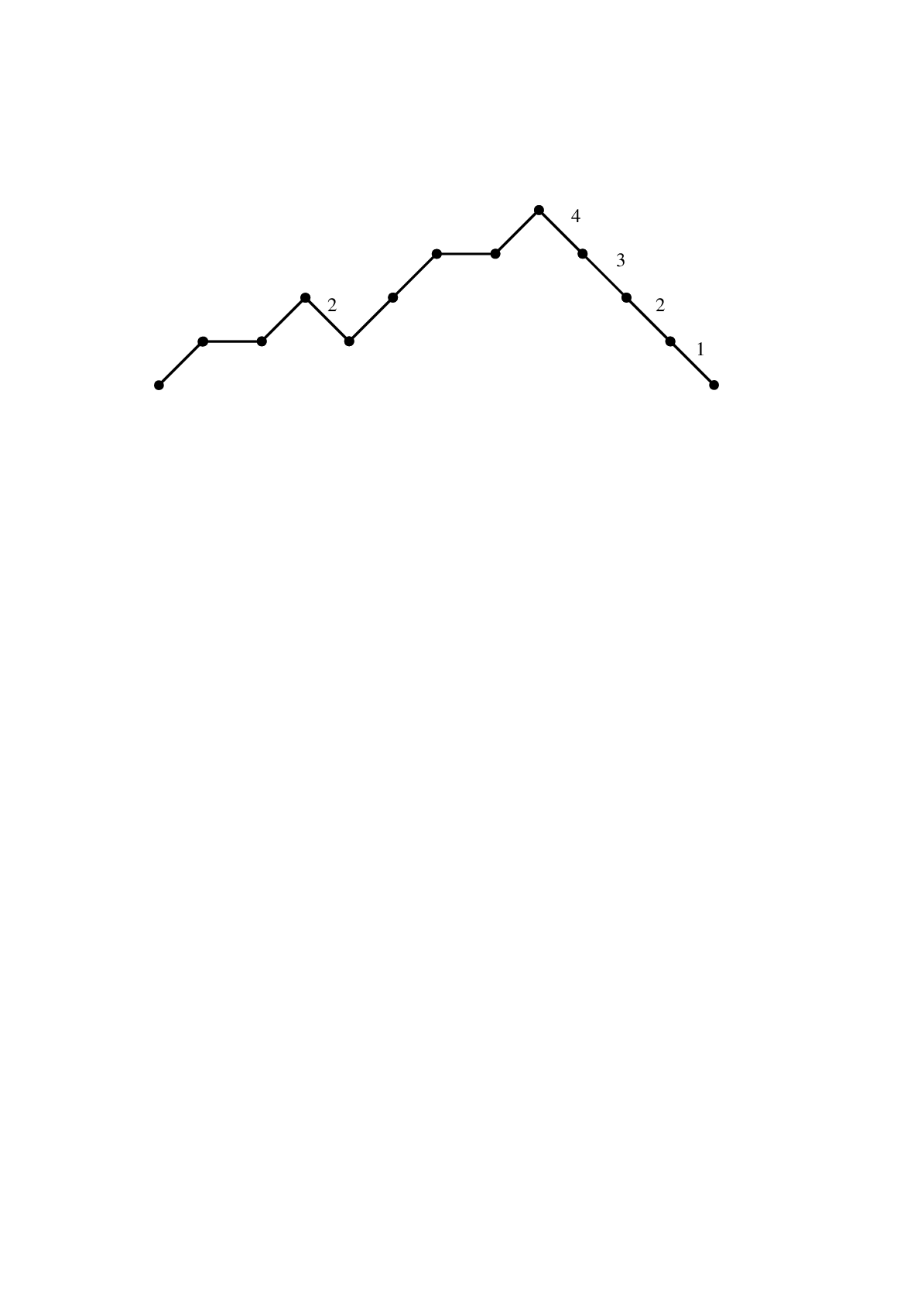}
\caption{The labelled Motzkin path associated with the $3412$
avoiding involution $\tau=12\ 2\ 4\ 3\ 11\ 10\ 7\ 9\ 8\ 6\ 5\ 1
$}\label{srof}
\end{center}
\end{figure}

\noindent The preceding result yields a bijection between the sets
$\mathscr{M}_n$ and $I_n(3412)$, that can be also found in
\cite{gui}.\newline

\noindent In \cite{egg} the cardinality of $I_n(3412,\tau)$ for
every $\tau$ in $S_3$ and $S_4$, and for several $\tau\in S_5$, is
determined. In the following, we exhibit a characterization for
labelled Motzkin paths corresponding to involutions in $I_n(3412)$
that avoid either a pattern $\pi\in S_3$ or $\pi=4321$. As a
consequence, we derive some of the enumerative results contained
in \cite{egg}.

\newtheorem{elpm}[yama]{Theorem}
\begin{elpm}\label{piude}
Let $\tau$ be an involution in $I_n(3412)$ associated with the
labelled Motzkin path $(M,\mu)$. Then:
\begin{itemize}
\item[i.] $\tau$ avoids $132$ if and only if $M$ does not contain either $HU$ or $DU$.
Similarly, $\tau$ avoids $213$ if and only if $M$ does not contain
either $DH$ or $DU$;
\item[ii.] $\tau$ avoids $321$ if and only if
  the irreducible components of $M$ are either $H$ or $UD$;
\item[iii.] $\tau$ avoids $312$, and hence $231$, if and only if
the irreducible components of $M$ are either of type $U^aHD^a$ or
$U^aD^a$;
\item[iv.] the set of paths corresponding to involutions in $I_n(3412)$ avoiding
$123$ can be constructed recursively as follows: either
$M=UH^aDUH^bD$, with $a+b=n-4$, or $M$ is obtained from a path of
the same kind by by prepending $U$ and appending $D$;
\item[v.] $\tau$ avoids $4321$ if and only if
  the irreducible components of $M$ are of type either $UH^aD$ or $H$.
\end{itemize}
\end{elpm}

\noindent \emph{Proof}
\begin{itemize}
\item[i.] Let $(M,\mu)$ be a Motzkin path
 with the maximal labelling. Then:
\begin{itemize}
\item[1.] if the $k$-th step in $M$ is horizontal and the
$(k+1)$-th step is an up step, then $\tau(k)=k$ and
$\tau(k+1)=j>k+1$. Hence, $\tau$ contains the $132$-pattern $k\ j\
k+1$;
\item[2.] if the $k$-th step in $M$ is a down step and the
$(k+1)$-th step is an up step, then $\tau(k)=j<k$ and
$\tau(k+1)=m>k+1$. Hence, $\tau$ contains the $132$-pattern $j\ m\
k+1$;
\end{itemize}
Suppose now that $M$ contains neither $HU$ nor $DU$. In this case,
$M=U^hS$, where $S$ is a suffix containing $h$ down steps and
$n-2h$ horizontal steps in some order. Then, the one-line notation
of $\tau$ has a prefix $e_1\cdots e_h$ consisting of all the
excedances in decreasing order and a suffix
$\tau(h+1)\cdots\tau(n)$ obtained by interlacing the sequence of
deficiences in decreasing order and the sequence of fixed points
in increasing order. It is easy to verify that such a $\tau$
avoids $132$, since each deficiency in $\tau$ must be less than
every fixed point and every excedance.

\begin{figure}[ht]
\begin{center}
\includegraphics[bb=61 563 470 655,width=.85\textwidth]{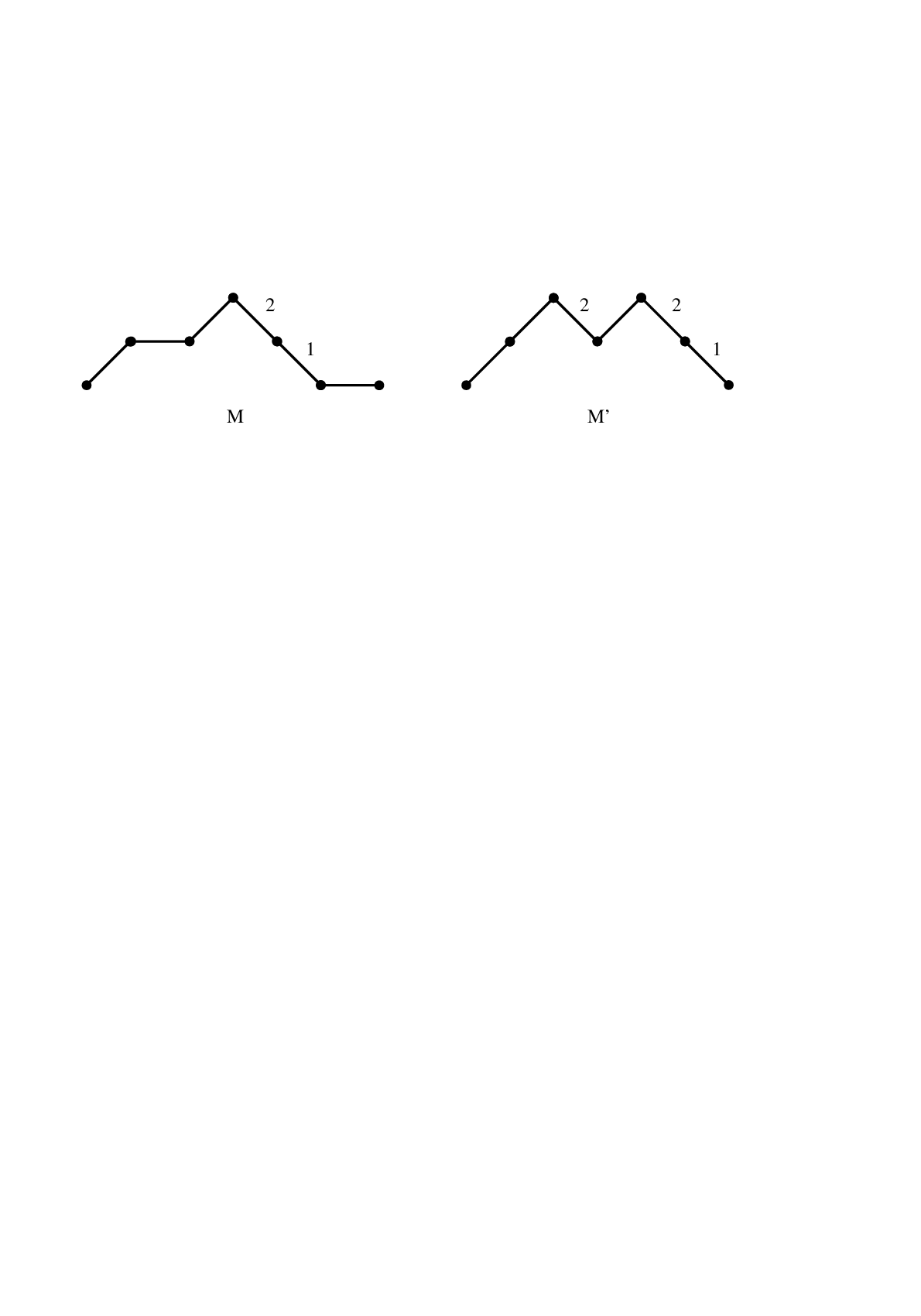} \caption{The labelled Motzkin paths $M$ and $M'$ corresponding to $\tau=5\ {\bf 2}\
{\bf 4}\ {\bf 3}\ 1\ 6$ and $\tau'=6\ 3\ {\bf 2}\ {\bf 5}\ {\bf 4
}\ 1$,
 that contain $132$.}
\end{center}
\end{figure}
\item[ii.] As stated in Proposition \ref{replace}, $321$-avoiding
involutions correspond to Motzkin paths whose irreducible
components are either Dyck paths or horizontal steps,  with the
unitary labelling. Hence, an involution in $I_n(3412,321)$ must
correspond to a Motzkin path of heigth at most $1$ whose
irreducible components are $H$ and $UD$, as desired.
\begin{figure}[ht]
\begin{center}
\includegraphics[bb=107 583 353 625,width=.5\textwidth]{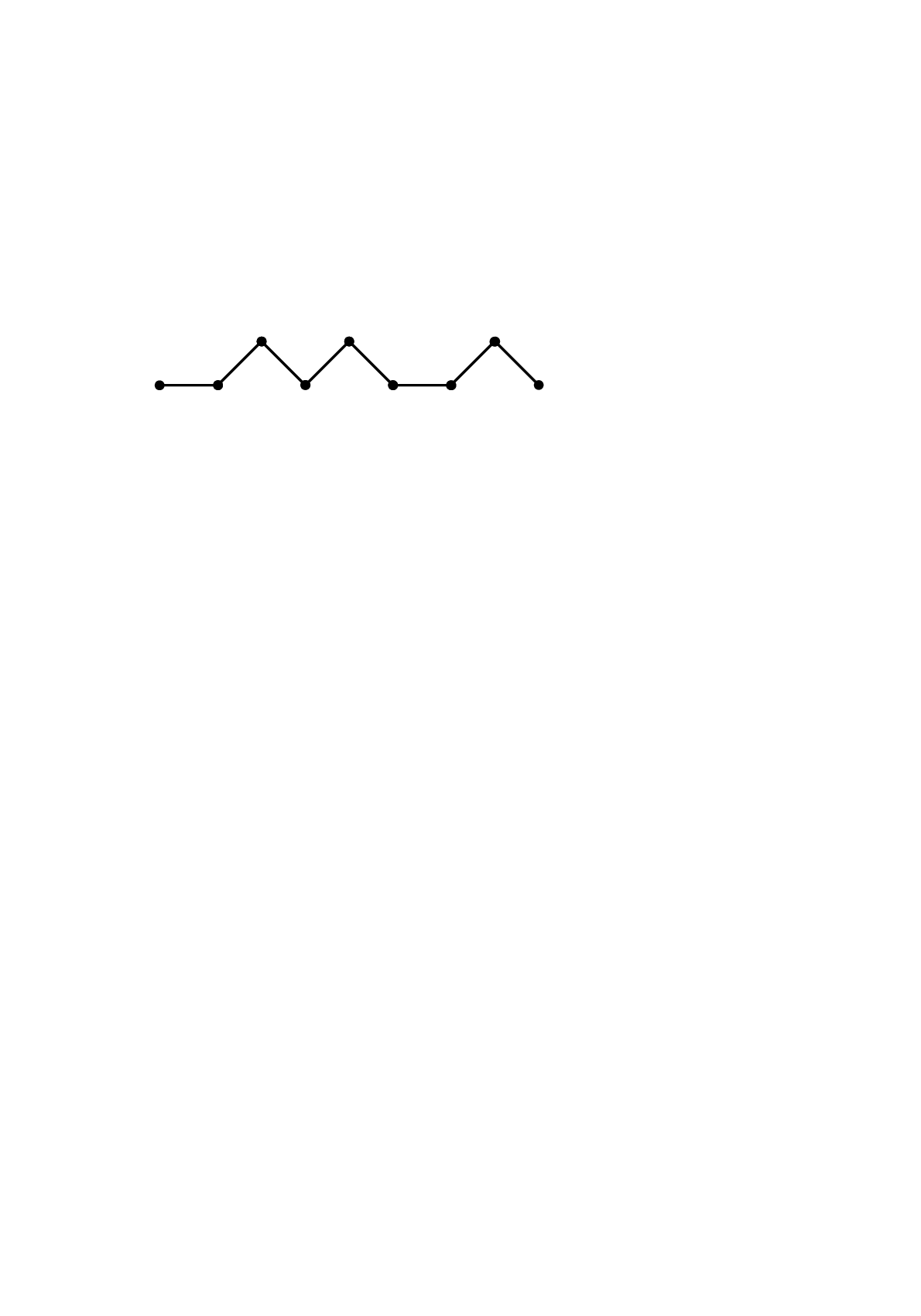} \caption{The labelled Motzkin path $M$ corresponding to $\tau=1\ 3\
2\ 5\ 4\ 6\ 8\ 7$,
 that avoids $321$.}
\end{center}
\end{figure}
\item[iii.] This characterization has been proved in Proposition \ref{fumetti}.
\item[iv.] Remark that an involution $\tau$ avoiding $123$ must have at
most $2$ connected components. If $\tau$ has exactly $2$
components, then there exists an integer $k$ such that the
one-line notation of $\tau$ is
$$k\ k-1\ \cdots\ 1\ n\ n-1\ \cdots\ k+1.$$
Hence, the labelled Motzkin path $\Phi(\tau)$ has two irreducible
components, of type either $U^aHD^a$ or $U^aD^a$,
 with the maximal labelling. It is easy to verify that
each one of these Motzkin paths corresponds to an involution
avoiding both $123$ and $3412$.

\noindent Consider now the case of a connected involution $\tau$
in $I_n(3412,123)$, corresponding to an irreducible Motzkin path
 with the maximal labelling. In this case, we must have
$\tau(1)=n$ and $\tau(n)=1$. Hence, the involution $\tau'\in
I_{n-2}$ obtained from $\tau$ by removing the symbols $1$ and $n$
and renormalizing the remaining symbols belongs to
$I_{n-2}(3412,123)$ if and only if $\tau$ belongs to
$I_n(3412,123)$.
\begin{figure}[ht]
\begin{center}
\includegraphics[bb=105 577 540 679,width=.85\textwidth]{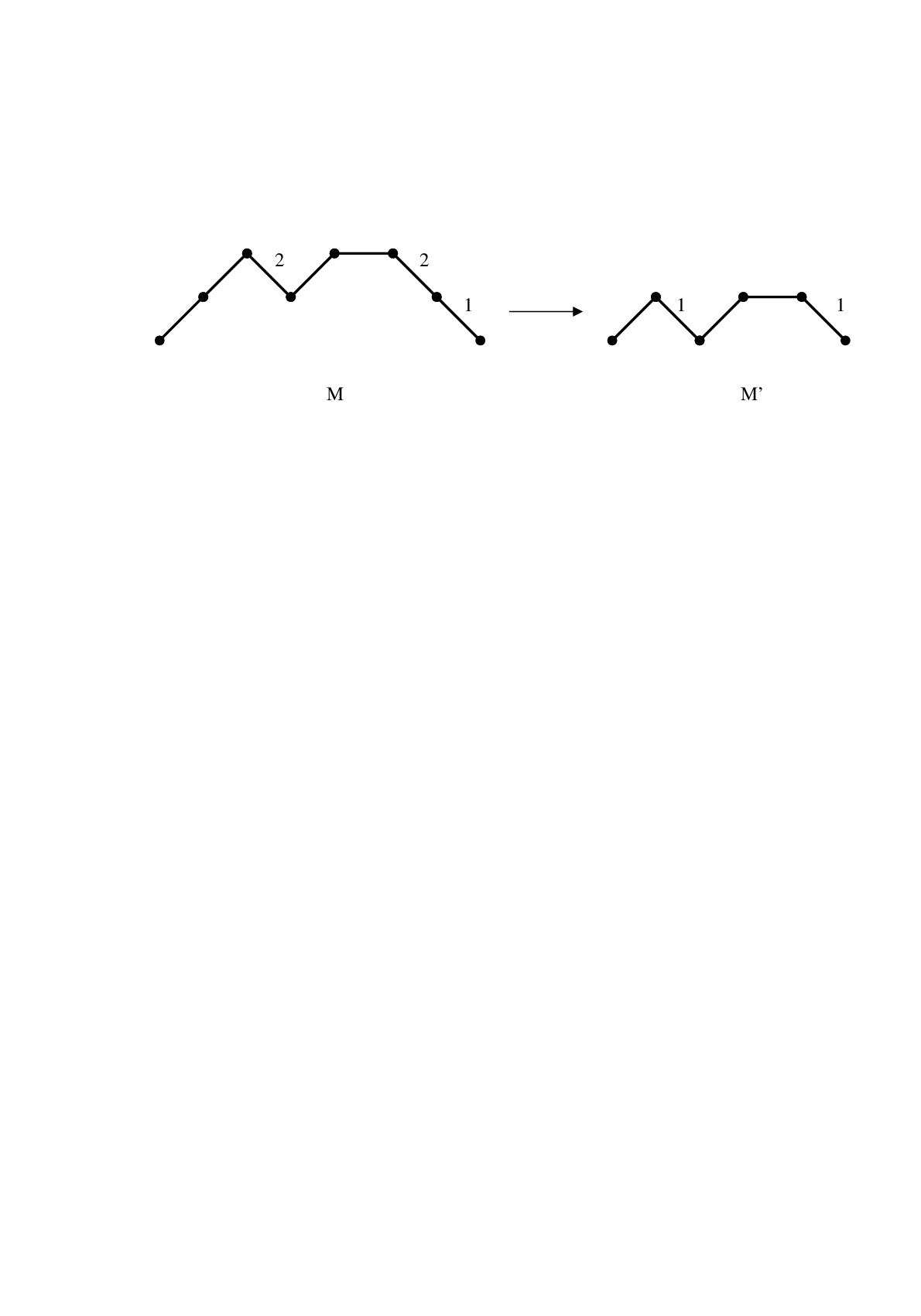} \caption{The Motzkin paths $M=\Phi(\tau)$ and $M'=\Phi(\tau')$,
where $\tau=7\ 3\ 2\ 6\ 5\ 4\ 1$
and $\tau'=2\ 1\ 5\ 4\ 3$ avoid $123$.}
\end{center}
\end{figure}
\item[v.] It is obvious that the height of a labelled Motzkin path corresponding
to an involution avoiding both $4321$ and $3412$ can not exceed
$1$. Then, its irreducible components must be of type either
$UH^aD$ or $H$.
\end{itemize}
\begin{flushright}
$\diamond$
\end{flushright}

\noindent The following enumerative results stated in \cite{egg}
can be now easily deduced from the preceding theorem:

\newtheorem{mple}[yama]{Corollary}
\begin{mple}\label{filip}
We have: \begin{itemize}
\item[i.]
$|I_n(3412,132)|=|I_n(3412,213)|=F_{n+1}$

\noindent where $F_n$ is the $n$-th Fibonacci number;
\item[ii.] $|I_n(3412,321)|=F_{n+1}$;
\item[iii.] $I_n(3412,312)=I_n(3412,231)=I_n(312)$ and
$|I_n(312)|=2^{n-1}$;
\item[iv.] $|I_n(3412,123)|=1+\lfloor\frac{n}{2}\rfloor\lceil\frac{n}{2}\rceil$
\item[v.] $|I_n(3412,4321)|=2^{n-1}$.
\end{itemize}
\end{mple}
\begin{flushright}
$\diamond$
\end{flushright}

\section{Fixed point free restricted involution}

\noindent Consider the set $DI_n$ of involutions without fixed
points on $n$ objects. It is evident that this set is nonempty if
and only if $n=2h$ is even. Involutions in $DI_n$ correspond, via
the bijection $\Phi$, to labelled Dyck paths.\newline

\noindent Theorem \ref{nlsn} allows us to
characterize involutions in $DI_n(4321)$ as follows:

\newtheorem{comemi}[yama]{Theorem}
\begin{comemi}\label{combi}
There is a bijection between the set $DI_{2h}(4321)$ and the set
of $\mathscr{D}_h$ of Dyck paths of length $2h$. Hence,
$|DI_{2h}(4321)|=C_h$, where $C_h$ is the $h$-th Catalan number.
\end{comemi}
\begin{flushright}
$\diamond$
\end{flushright}

\noindent It has been proved (see \cite{drs}) that also the
cardinality of the set $DI_{2h}(321)$ is $C_h$. This fact is
somehow surprising since $DI_n(321)\subseteq DI_n(4321)$. The
described bijection yields an easy proof of the fact that these
two sets coincide.

\newtheorem{infa}[yama]{Proposition}
\begin{infa}
We have:$$DI_{2h}(321)=DI_{2h}(4321).$$
\end{infa}

\noindent \emph{Proof} The characterization of labelled Motzkin
paths associated with involutions in $I_{2h}(321)$ given in
Proposition \ref{replace} implies that the set $DI_{2h}(321)$
corresponds bijectively to the set of Dyck paths of semilength $h$
 with the unitary labelling. The assertion now follows from
Theorem \ref{combi}.
\begin{flushright}
$\diamond$
\end{flushright}

\noindent The study of fixed point free involutions that avoid
$4321$ and a further pattern $\pi\in S_3$, with $\pi\neq 321$, has
no interest. In fact, Theorem \ref{attrezz} implies that
$|DI_n(4321,\pi)|\leq 2$ for every $n\in\mathbb{N}$.\newline

\noindent Consider now the set $DI_{2h}(3412)$ of fixed point free
involutions avoiding $3412$. Also in this case, by Theorem
\ref{consimm}, the set $DI_{2h}(3412)$ is in bijection with the
set of Dyck paths of lenght $2h$, whose cardinality is $C_h$.

\noindent The characterization of the sets $DI_{2h}(3412,\pi)$,
where $\pi$ is one of the patterns considered in Theorem
\ref{piude}, is nontrivial only in the following cases:

\newtheorem{following}[yama]{Theorem}
\begin{following}
For every integer $n$, we have:
\begin{itemize}
\item[a.] $|DI_{2h}(3412,123)|=1+{h \choose 2}$;
\item[b.] $|DI_{2h}(3412,312)|=2^{h-1}$.
\end{itemize}
\end{following}

\noindent \emph{Proof}
\begin{itemize}
\item[a.] Non-connected involutions in $DI_{2h}(3412,123)$ correspond bijectively to Dyck paths of type
$D=U^aD^aU^{h-a}D^{h-a}$, where $0<a<h$. There are $h-1$
involutions of this type. As seen in the proof of Theorem
\ref{piude}.iv, connected involutions in $DI_{2h}(3412,123)$ are
in bijection with involutions in $DI_{2h-2}(3412,123)$. Hence, if
we set $d_{2h}=|DI_{2h}(3412,123)|$, we have:
$$d_{2h}=d_{2h-2}+h-1.$$
Since $d_2=1$, we get the assertion.
\item[b.] Remark that $DI_{2h}(3412,312)=DI_{2h}(312)$.
Involutions in $I_{2h}(312)$ without fixed points correspond
bijectively to Dyck paths whose irreducible components are of type
$U^aD^a$ (see Theorem \ref{fumetti}). These paths are encoded by
the compositions of $2h$ into even parts.
\end{itemize}
\begin{flushright}
$\diamond$
\end{flushright}

\section{Restricted centrosymmetric involutions}

\noindent A \emph{centrosymmetric involution} is an element
$\tau\in I_n$ such that $\tau_{rc}=\tau$, namely, an involution
such that $\tau(i)+\tau(n+1-i)=n+1$, for every $1\leq i\leq n$.
Denote by $CI_n$ the set of
centrosymmetric involutions in $I_n$.\\

\noindent The study of pattern avoidance on the set of
centrosymmetric involutions has been carried out by Egge
\cite{eggedue} for every pattern $\pi\in S_3$. Moreover, Guibert
and Pergola considered in \cite{pg} the case of vexillary
centrosymmetric involutions, namely,
centrosymmetric involutions avoiding $2143$. In this section we specialize the previous results to the set $CI_n$.\\

\noindent In the case of centrosymmetric involutions, we have
further equidistribution results among patterns. In fact, $CI_n$
is closed under both reverse and complement operations, that was
not the case for the set $I_n$. This implies that the patterns
$\pi$, $\pi_c$,
$\pi_r$, and $\pi_{rc}$ are equidistributed over $CI_n$.\\

\noindent Theorem \ref{modo} implies that the Motzkin path
associated with a centrosymmetric involution $\tau\in CI_n$ must
be symmetric with respect to the vertical line $x=\frac{n}{2}$.

\begin{figure}[ht]
\begin{center}
\includegraphics[bb=57 611 377 704,width=.6\textwidth]{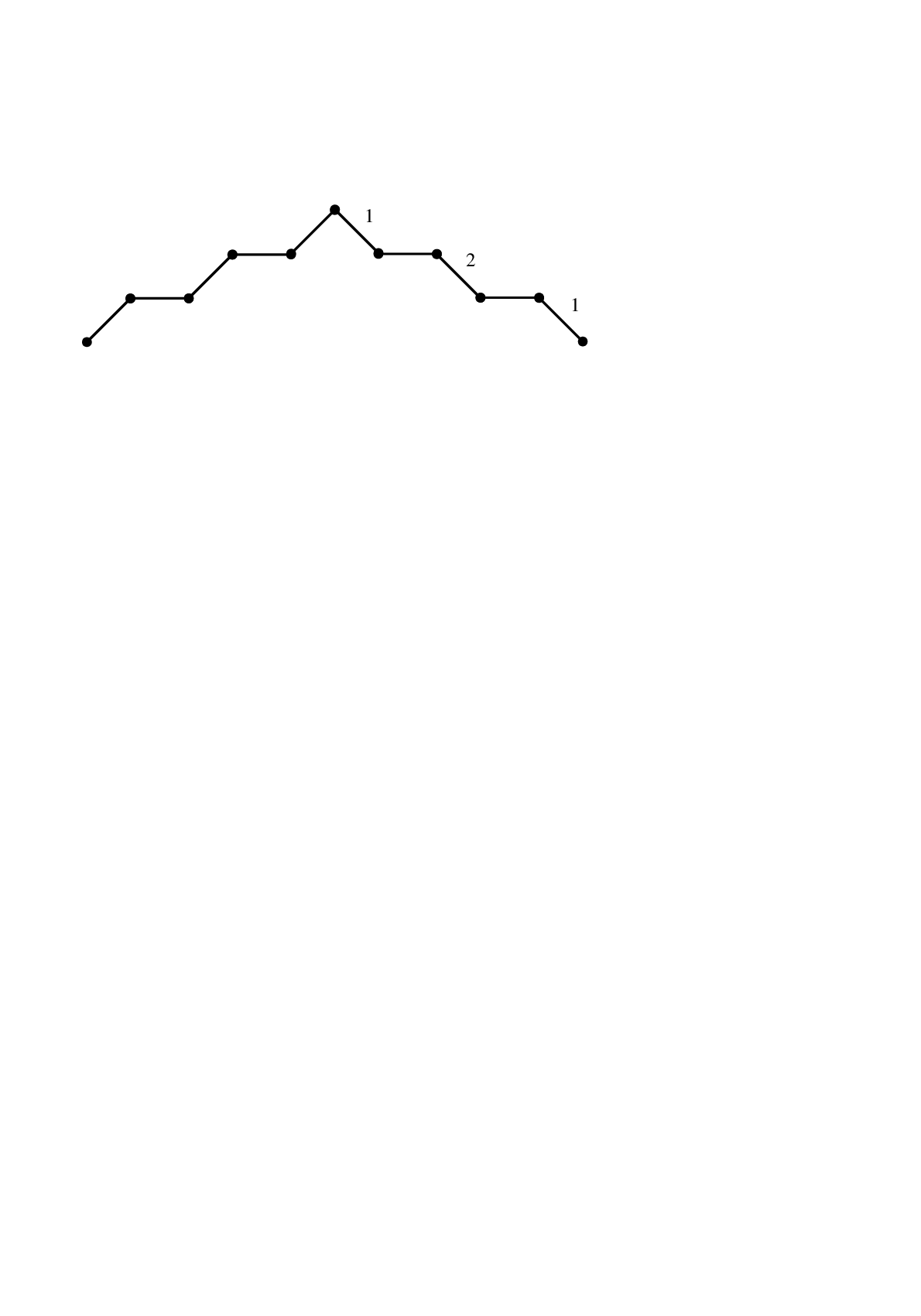}
\caption{A symmetric labelled Motzkin path that is not associated with a centrosymmetric
involution.}\label{inve}
\end{center}
\end{figure}

\noindent The converse is false, in general. In fact, for
instance, the involution associated with the symmetric labelled
Motzkin path in Figure \ref{inve} is $$\tau=6\ 2\ 10\ 4\ 8\ 1\ 7\
5\ 9\ 3,$$ that is not centrosymmetric.

\noindent However, if a symmetric Motzkin path $M$ is endowed with
either the maximal or the unitary labelling, then the
corresponding involution is centrosymmetric. In fact, consider for
example the case of the unitary labelling, and suppose that there
exists a symmetric Motzkin path $M$ such that the corresponding
involution $\tau\in I_n(4321)$ is not centrosymmetric, namely,
$\tau_{rc}\neq \tau$. Then, also $\tau_{rc}$ avoids $4321$, as
seen in the proof of Theorem \ref{attrezz}.i. This implies that
$\tau_{rc}$ corresponds to the same Motzkin path with the same
labelling, hence contradicting Theorem \ref{nlsn}.

\noindent The case of the maximal labelling can be treated
similarly, remarking that also the set $I_n(3412)$ is closed under
the
reverse-complement map.\\

\noindent The bijection between Motzkin paths and centrosymmetric
involutions avoiding $4321$ (or $3412$) yields the following
enumerative results:

\newtheorem{unaltro}[yama]{Theorem}
\begin{unaltro}
For every integer $h$,
\begin{equation}|CI_{2h+1}(4321)|=|CI_{2h+1}(3412)|=|CI_{2h}(4321)|=$$$$=|CI_{2h}(3412)|=\sum_{i=0}^h
\frac{h!}{(h-i)!\left\lfloor\frac{i}{2}\right\rfloor!\left\lceil\frac{i}{2}\right\rceil!}.\label{hairag}\end{equation}
\end{unaltro}

\noindent \emph{Proof} First of all, remark that a symmetric
Motzkin path of length $2h+1$ can be uniquely obtained by adding a
horizontal step in the middle position of a symmetric Motzkin path
of length $2h$. Hence, we can restrict our attention to the even
case. A symmetric Motzkin path of length $2h$ is completely
determined by its first $h$ steps. This means that we only need to
count Motzkin prefixes of length $h$, whose formula can be found,
for example, in \cite{pg}.
\begin{flushright}
$\diamond$
\end{flushright}

\noindent The following double restriction results can be derived
from the analogous theorems for the set of all involutions:

\newtheorem{rifiu}[yama]{Theorem}
\begin{rifiu}
We have:
\begin{itemize}
\item[i.] $|CI_n(4321,132)|=\left\lfloor\frac{n}{2}\right\rfloor$;
\item[ii.] $CI_n(4321,321)=CI_n(321)$. Moreover, we have: \begin{itemize}\item[]$|CI_{2h}(4321,321)|=2^h$ \item[]$|CI_{2h+1}(4321,321)|={h\choose
\left\lfloor\frac{h}{2}\right\rfloor}$;\end{itemize}
\item[iii.] $|CI_n(4321,312)|=t_{\left\lfloor\frac{n}{2}\right\rfloor+1}+t_{\left\lfloor\frac{n}{2}\right\rfloor+2}$,\\ where $t_k$
denotes the $k$-th tribonacci number;
\item[iv.]
$|CI_n(3412,132)|=\left\lfloor\frac{n}{2}\right\rfloor$;
\item[v.]
$|CI_{2h}(3412,321)|=F_{h+2}$ \hspace{2cm}
$|CI_{2h+1}(3412,321)|=F_{h+1}$;
\item[vi.] $CI_n(3412,312)=CI_n(312)$ \hspace{2cm}
$|CI_n(3412,312)|=2^{\left\lfloor\frac{n}{2}\right\rfloor}$;
\item[vii.] $|CI_{2h}(3412,123)|=h+1$ \hspace{1.7cm}
$|CI_{2h+1}(3412,123)|=1$;
\item[viii.]
$|CI_n(3412,4321)|=2^{\left\lfloor\frac{n}{2}\right\rfloor}$.
\end{itemize}
\end{rifiu}

\noindent \emph{Proof}

\begin{itemize}
\item[i.] as stated in Theorem \ref{attrezz}.i, involutions avoiding both $4321$ and $132$ are in bijection with
Motzkin paths of the kind $M=U^aH^bD^aH^c$. Among these, only the
$\left\lfloor\frac{n}{2}\right\rfloor$ with $c=0$ are symmetric.
\item[ii.] See \cite{eggedue}.
\item[iii.] Involutions in $CI_n(4321,312)$ correspond bijectively to symmetric \linebreak Motzkin
paths whose irreducible components are either $H$ or $UD$ or
$UHD$. When $n$ is odd, a path of this kind can be uniquely
constructed by adding a horizontal step at the middle position of
a Motkin path of the same kind of length $n-1$. Hence, we consider
only the even case $n=2h$. The first half of a Motzkin path of
this type is either a (possibly non symmetric) Motzkin path of
this type of length $h$ or a Motzkin path of the prescribed type
of length $h-1$ followed by an up step. This completes the proof.
\item[iv.] In this case, by Thoerem \ref{piude}.i, involutions avoiding both $3412$ and $132$ are in bijection with
Motzkin paths of the kind $M=U^aS$. Among them, only the
$\left\lfloor\frac{n}{2}\right\rfloor$  paths of type
$M=U^aH^bD^a$ are symmetric.
\item[v.] Involutions in $CI_n(3412,321)$ correspond bijectively to symmetric \linebreak Motzkin paths
whose irreducible components are either $H$ or $UD$. Consider first the even case
$n=2h$. The left half of a symmetric Motzkin path of the
prescribed type can be either a (possibly non symmetric) Motzkin
path of the same type of length $h$, or a Motzkin path of the same
type of length $h-1$, followed by an up step $U$. Hence, the
number of these paths is $F_{h+1}+F_h=F_{h+2}$, by Theorem
\ref{filip}.ii.

\noindent On the other hand, if $n=2h+1$, Motzkin paths of the
prescribed type must be decomposable as follows
$$M=NHN,$$
where $N$ is a Motzkin path whose connected components are either
$H$ or $UD$. Hence,  by Theorem \ref{filip}.ii, the number of
these paths is $F_{h+1}$.
\item[vi.] See \cite{eggedue}.
\item[vii.] Consider first the odd case. It is easily seen that all
involutions in \linebreak $CI_{2h+1}(3412,123)$ are connected.
Hence, as seen in the proof of Theorem \ref{filip}.iv, a Motzkin
path associated with an involution in \linebreak
$CI_{2h+1}(3412,123)$ can be uniquely obtained by prepending $U$
and appending $D$ to a Motzkin path associated with an involution
in \linebreak $CI_{2h-1}(3412,123)$. Since $CI_{1}(3412,123)$
contains only one element, we have $|CI_{2h+1}(3412,123)|=1$.

\noindent If $n=2h$, all involutions in $CI_{2h}(3412,123)$ are
connected, except for the involution corresponding to either
$U^aHD^aU^aHD^a$ or $U^aD^aU^aD^a$, where
$a=\left\lfloor\frac{h}{2}\right\rfloor$. Connected involutions in
$CI_{2h}(3412,123)$ correspond bijectively to the elements of
$CI_{2h-2}(3412,123)$, as remarked above. This implies that,
setting $c_h=|CI_{2h}(3412,123)|$, we have $c_h=c_{h-1}+1$. Since
$c_1=1$, we get the assertion.
\item[viii.] The elements in $CI_n(3412,4321)$ correspond to symmetric Motzkin paths of
height at most $1$. A path $M$ of this type is completely
determined by the prefix $p(M)$ consisting of the
$\left\lfloor\frac{n}{2}\right\rfloor$ leftmost steps of $M$. Fix
a non negative integer $h\leq
\left\lfloor\frac{n}{2}\right\rfloor$. Motzkin prefixes whose
rightmost point on the $x$-axis is $(h,0)$ correspond bijectively
to Motzkin paths of length $h$ and height $1$. In fact, such a
prefix uniquely decomposes into $p(M)=M'UH^{b}$, where
$b=\left\lfloor\frac{n}{2}\right\rfloor-h-1$. Hence, we have:
$$|CI_n(3412,4321)|=1+\sum_{h=1}^{\left\lfloor\frac{n}{2}\right\rfloor} 2^{h-1}=2^{\left\lfloor\frac{n}{2}\right\rfloor}.$$
\end{itemize}

\begin{flushright}
$\diamond$
\end{flushright}

\end{document}